\documentclass[12pt]{article}
\usepackage{amsmath,amsfonts,amssymb,amsthm}
\usepackage{color}
\usepackage{colortbl}
\usepackage{array}
\usepackage{mathrsfs}
\usepackage{dcolumn}
\usepackage{longtable}
\usepackage{hhline}
\usepackage{bm}

\newcommand{\braced}[2]{\genfrac{\{}{\}}{0pt}{0}{#1}{#2}}

\theoremstyle{plain}
\newtheorem{thm}{Theorem}[section]

\newtheorem{cor}[thm]{Corollary}

\theoremstyle{definition}
\newtheorem{defn}[thm]{Definition}
\newtheorem{rem}[thm]{Remark}

\title{\bf On Multi Poly-Bernoulli Polynomials}
\author{
{\normalsize Roberto B. Corcino, Hassan Jolany, Cristina B. Corcino and Takao Komatsu
}
}

\date{}

\begin{document}

\maketitle

\begin{abstract}
In this paper, we define multi poly-Bernoulli polynomials using multiple polylogarithm and derive some properties parallel to those of poly-Bernoulli polynomials. Furthermore, an explicit formula for certain Hurwitz-Lerch type multi poly-Bernoulli polynomials is established using the $r$-Whitney numbers of the second kind.    


\end{abstract}

\section{Introduction}

The poly-Bernoulli numbers were first introduced by Kaneko \cite{Kaneko} by means of the following exponential generating function
$$\frac{{\rm Li}_k(1-e^{-t})}{1-e^{-t}}=\sum_{n=0}^{\infty}B^{(k)}_n\frac{t^n}{n!}$$  
where
$${\rm Li}_k(z)=\sum_{n=1}^{\infty}\frac{z^n}{n^k}.$$
One can easily verify that, when $k=1$, this gives $B_n^{(1)}=B_n(1)$, a particular value of Bernoulli polynomial. The poly-Bernoulli numbers can be extended in polynomial form as
$$\frac{{\rm Li}_k(1-e^{-t})}{1-e^{-t}}e^{xt}=\sum_{n=0}^{\infty}B^{(k)}_n(x)\frac{t^n}{n!}.$$
Several properties for these numbers and polynomials are established including the following explicit formula of Arakawa and Kaneko \cite{Arakawa} 
\[B^{(-k)}_n=\sum_{m\geq0}m!\braced{n+1}{m+1}m!\braced{k+1}{m+1},
\]
where $\braced{n}{m}$ denote Stirling numbers of the second kind. Through this formula, one may be able to interpret the numbers $B^{(-k)}_n$ in terms of the number of partitions as $\braced{n}{m}$ count the number of partitions of an $n$-set into $m$ nonempty subsets. Recently, these numbers have been interpreted as the number of binary lonesum matrices of size $n\times k$ where the binary lonesum matrix is a binary matrix which can be reconstructed from its row and column sums \cite{Benyi}.

\smallskip
It is worth-mentioning that the above generalization of Kaneko has been generalized further by Cenkci and Young \cite{Cenkci} using the concept of Hurwitz-Lerch zeta function $\Phi (z,s,a)$ as follows
\begin{equation}\label{HpolyB}
\Phi (1-e^{-t},k,a)=\sum_{n=0}^{\infty}B_{n,a}^{(k)}\frac{t^n}{n!}
\end{equation}
where
\begin{equation}\label{HLerch}
\Phi (z,k,a)=\sum_{n=0}^{\infty}\frac{z^n}{(n+a)^k}.
\end{equation}
The numbers $B_{n,a}^{(k)}$ are called the Hurwitz type poly Bernoulli numbers. These numbers have been shown Theorem 2.1 of \cite{Cenkci} to have explicit formula
\begin{equation}\label{HLExplicit}
B_{n,a}^{(k)}=(-1)^n\sum_{m=0}^{n}\frac{(-1)^mm!\braced{n}{m}}{(m+a)^k}.
\end{equation}
Obviously, $B_{n,a}^{(-k)}$ are nonnegative integers when $a$ is a nonnegative integer. Hence, combinatorial interpretation for these numbers is possible to establish.  

\smallskip
Recently, several mathematicians and physicists have been attracted to work on nested harmonic sums because of their rich structure that fascinates theoreticians working on the areas like algebra, number theory and combinatorics. These sums occur in the study on knot theory and quantum field theory. Multiple polylogarithms are certain generalization of the nested harmonic sums as well as the Riemann zeta function and the ordinary polylogarithm, which preserve many interesting properties. More precisely, the multiple polylogarithms are defined by
\begin{equation}\label{multipolylogNew}
{\rm Li}_{(k_1,k_2,\ldots, k_r)}(z_1, z_2, \ldots, z_r)=\sum_{ 0< m_1< m_2<\ldots < m_r }\prod_{j=1}^rm_j^{-k_j}z_j^{n_j},
\end{equation}
where $k_1$, $k_2$, $\ldots$, $k_r$ and $z_1$, $z_2$, $\ldots$, $z_r$ are complex numbers suitably restricted so that the
sum (\ref{multipolylogNew}) converges. These polynomials also occur in various fields like combinatorics, knot theory, quantum field theory and mirror symmetry. In fact, there are several sophisticated study that relate multiple polylogarithms to arithmetic and algebraic geometry and to algebraic $K$-theory.

\smallskip
There are studies that relate a special case of the aforementioned multiple polylogarithms to the concept of Bernoulli and Euler numbers and polynomials by extending the idea of poly-Bernoulli and poly-Euler numbers and polynomials to multiple parameter case. In particular, some consider the case in which 
$$z_1 = z_3 = \ldots = z_{r-1} = 1, z_r=z,$$
that is, the studies that deal with the multiple polylogarithm of the form
\begin{equation}\label{multipolylog}
{\rm Li}_{(k_1,k_2,\ldots, k_r)}(z)=\sum_{ 0< m_1< m_2<\ldots < m_r }\frac{z^{m_r}}{m_1^{k_1} m_2^{k_2}\ldots m_r^{k_r}}.
\end{equation}
For instance, Imatomi et al. \cite{Imatomi} have defined a certain generalization of Bernoulli numbers in terms of these multiple logarithms as follows 
\begin{equation}\label{MultiBernoulli}
\frac{{\rm Li}_{(k_1, k_2,\ldots, k_r)}(1-e^{-t})}{1-e^{-t}}=\sum_{n=0}^{\infty}{B}^{(k_1, k_2,\ldots, k_r)}_n\frac{t^n}{n!}.
\end{equation}
These numbers possess respectively the following recurrence relation and explicit formula
\begin{align}
{B}^{(k_1, k_2,\ldots, k_r)}_n&=\frac{1}{n+1}\left({B}^{(k_1-1, k_2,\ldots, k_r)}_n-\sum_{m=1}^{n-1}\binom{n}{m-1}{B}^{(k_1, k_2,\ldots, k_r)}_m\right)\\
{B}^{(k_1, k_2,\ldots, k_r)}_n&=(-1)^n\sum_{n+1\geq m_1>m_2>\ldots >m_r>0}\frac{(-1)^{m_1-1}(m_1-1)!S(n,m_1-1)}{m_1^{k_1}m_1^{k_2}\ldots m_r^{k_r}}.
\end{align}
Parallel to the above generalization is the generalized multi poly-Euler polynomials which are denoted by ${E}^{(k_1, k_2,\ldots, k_r)}_n(x;a,b,c)$. These polynomials have been introduced in \cite{CorJol1} by means of the above multiple poly-logarithm, also known as multiple zeta values. More precisely, we have
\begin{equation}\label{multipolyeuler}
\frac{2{\rm Li}_{(k_1, k_2,\ldots, k_r)}(1-(ab)^{-t})}{(a^{-t}+b^{t})^r}c^{rxt}=\sum_{n=0}^{\infty}{E}^{(k_1, k_2,\ldots, k_r)}_n(x;a,b,c)\frac{t^n}{n!}.
\end{equation}
When $r=1$, (\ref{multipolyeuler}) boils down to the generalized poly-Euler polynomials with three parameters $a, b, c$. Moreover, when $c=e$, (\ref{multipolyeuler}) reduces to the multi poly-Euler polynomials with two parameters $a, b$. These special cases have been discussed intensively in \cite{CorJol1, CorJol2}.

\smallskip
In this paper, generalized multi poly-Bernoulli polynomials are defined and some properties of these polynomials are established parallel to those of the poly-Bernoulli polynomials. Moreover, certain generalization of multi poly-Bernoulli numbers is defined in terms of generalized Hurwitz-Lerch multiple zeta values.

\section{Generalized Multi Poly-Bernoulli Polynomials}
Parallel to the definition of generalized multi poly-Euler polynomials in (\ref{multipolyeuler}), we have the following generalization of poly-Bernoulli numbers.

\begin{defn}\label{multipolybernoulli}\rm
The generalized multi poly-Bernoulli polynomials are defined by
\begin{equation}\label{multipolybern}
\frac{{\rm Li}_{(k_1, k_2,\ldots, k_r)}(1-(ab)^{-t})}{(b^{t}-a^{-t})^r}c^{rxt}=\sum_{n=0}^{\infty}{B}^{(k_1, k_2,\ldots, k_r)}_n(x;a,b,c)\frac{t^n}{n!}.
\end{equation}
\end{defn}
The following theorem contains some identities for ${B}^{(k_1, k_2,\ldots, k_r)}_n(x;a,b,c)$.

\begin{thm}\label{addd} The generalized multi poly-Bernoulli polynomials satisfy the following identities.
\begin{align*}
{B}^{(k_1, k_2,\ldots, k_r)}_n(x;a,b,c)&=\sum_{m=0}^{\infty}\sum_{l=m}^n(r\log c)^l\braced{l}{m}\binom{n}{l}{B}^{(k_1, k_2,\ldots, k_r)}_{n-l}(-m\log c;a,b)(x)^{(m)}\\
{B}^{(k_1, k_2,\ldots, k_r)}_n(x;a,b,c)&=\sum_{m=0}^{\infty}\sum_{l=m}^n(r\log c)^l\braced{l}{m}\binom{n}{l}{B}^{(k_1, k_2,\ldots, k_r)}_{n-l}(0;a,b)(x)^{(m)}\\
{B}^{(k_1, k_2,\ldots, k_r)}_n(x;a,b,c)&=\sum_{m=0}^{\infty}\binom{n}{m}\sum_{l=0}^{n-m}\frac{\binom{n-m}{l}}{\binom{l+s}{l}}\braced{l+s}{s}{B}^{(k_1, k_2,\ldots, k_r)}_{n-m-l}(0;a,b)B^{(s)}_m(xr\log c)\\
{B}^{(k_1, k_2,\ldots, k_r)}_n(x;a,b,c)&=\sum_{m=0}^{n}\frac{\binom{n}{m}}{(1-\lambda)^s}\sum_{j=0}^s\binom{s}{j}(-\lambda)^{s-j}{B}^{(k_1, k_2,\ldots, k_r)}_{n-m}(j;a,b)H^{(s)}_m(xr\log c;\lambda)
\end{align*}
where $(x)^{(n)}=x(x+1)\ldots (x+n-1)$, $(x)_n=x(x-1)\ldots (x-n+1)$,
$$\left(\frac{t}{e^t-1}\right)^se^{xt}=\sum_{n=0}^{\infty}B^{(s)}_n(x)\frac{t^n}{n!}\mbox{  and  }\left(\frac{1-\lambda}{e^t-\lambda}\right)^se^{xt}=\sum_{n=0}^{\infty}H^{(s)}_n(x;\lambda)\frac{t^n}{n!}.$$
\begin{proof}
For the first identity, note that (\ref{multipolybern}) can be written as
$$\sum_{n=0}^{\infty}{B}^{(k_1, k_2,\ldots, k_r)}_n(x;a,b,c)\frac{t^n}{n!}=\frac{{\rm Li}_{(k_1, k_2,\ldots, k_r)}(1-(ab)^{-t})}{(b^{t}-a^{-t})^r}(1-(1-e^{-rt\log c}))^{-x}.$$
Using Newton's Binomial Theorem, we have
$$\sum_{n=0}^{\infty}{B}^{(k_1, k_2,\ldots, k_r)}_n(x;a,b,c)\frac{t^n}{n!}=\frac{{\rm Li}_{(k_1, k_2,\ldots, k_r)}(1-(ab)^{-t})}{(b^{t}-a^{-t})^r}\sum_{m=0}^{\infty}\binom{x+m-1}{m}(1-e^{-rt\log c})^m$$
\begin{align*}
&=\sum_{m=0}^{\infty}(x)^{(m)}\frac{(e^{rt\log c}-1)^m}{m!}\frac{{\rm Li}_{(k_1, k_2,\ldots, k_r)}(1-(ab)^{-t})}{(b^{t}-a^{-t})^r}e^{-mrt\log c}\\
&=\sum_{m=0}^{\infty}(x)^{(m)}\left(\sum_{n=0}^{\infty}\braced{n}{m}\frac{(rt\log c)^n}{n!}\right)\left(\sum_{n=0}^{\infty}{B}^{(k_1, k_2,\ldots, k_r)}_n(-mr\log c;a,b)\frac{t^n}{n!}\right)\\
&=\sum_{n=0}^{\infty}\left\{\sum_{m=0}^{\infty}\sum_{l=m}^{n}(r\log c)^l\braced{l}{m}\binom{n}{l}{B}^{(k_1, k_2,\ldots, k_r)}_{n-l}(-mr\log c;a,b)(x)^{(m)}\right\}\frac{t^n}{n!}.
\end{align*}
Comparing coefficients completes the proof of the first identity. For the next identity, this can be shown parallel to the above argument. For the last two identities, note that (\ref{multipolybern}) can be written as
$$G(t)=\sum_{n=0}^{\infty}{B}^{(k_1, k_2,\ldots, k_r)}_n(x;a,b,c)\frac{t^n}{n!}\qquad\qquad\qquad\qquad\qquad\qquad\qquad\qquad\qquad\qquad\qquad\qquad\qquad\qquad\qquad\qquad$$
\begin{align*}
&=\left(\frac{(e^t-1)^s}{s!}\right)\left(\frac{t^se^{xrt\log c}}{(e^t-1)^s}\right)\left(\frac{2{\rm Li}_{(k_1, k_2,\ldots, k_r)}(1-(ab)^{-t})}{(b^{t}-a^{-t})^r}\right)\frac{s!}{t^s}\\
&=\left(\sum_{n=0}^{\infty}\braced{n+s}{s}\frac{t^{n+s}}{(n+s)!}\right)\left(\sum_{m=0}^{\infty}B^{(s)}_m(xr\log c)\frac{t^m}{m!}\right)\left(\sum_{n=0}^{\infty}{B}^{(k_1, k_2,\ldots, k_r)}_n(0;a,b)\frac{t^n}{n!}\right)\frac{s!}{t^s}\\
&=\sum_{n=0}^{\infty}\left\{\sum_{m=0}^{n}\binom{n}{m}\sum_{l=0}^{n-m}\frac{\binom{n-m}{l}}{\binom{l+s}{s!}}\braced{l+s}{s}{B}^{(k_1, k_2,\ldots, k_r)}_{n-m-l}(0;a,b)B^{(s)}_m(xr\log c)\right\}\frac{t^n}{n!}.
\end{align*}
This proves the third identity. Again, using (\ref{multipolybern}), we have
$$\sum_{n=0}^{\infty}{B}^{(k_1, k_2,\ldots, k_r)}_n(x;a,b,c)\frac{t^n}{n!}=\left(\frac{(1-\lambda)^s}{(e^t-\lambda)^s}e^{xrt\log c}\right)\left(\frac{(e^t-\lambda)^s}{(1-\lambda)^s}\right)\left(\frac{2{\rm Li}_{(k_1, k_2,\ldots, k_r)}(1-(ab)^{-t})}{(b^{t}-a^{-t})^r}\right)$$
\begin{align*}
&=\left(\sum_{n=0}^{\infty}H^{(s)}_n(xr\log c;\lambda)\frac{t^n}{n!}\right)\left(\sum_{j=0}^{s}\binom{s}{j}(-\lambda)^{s-j}\frac{2{\rm Li}_{(k_1, k_2,\ldots, k_r)}(1-(ab)^{-t})}{(b^{t}-a^{-t})^r}e^{jt}\right)\\
&=\sum_{n=0}^{\infty}\left(\sum_{m=0}^{n}\frac{\binom{n}{m}}{(1-\lambda)^s}\sum_{j=0}^{s}\binom{s}{j}(-\lambda)^{s-j}{B}^{(k_1, k_2,\ldots, k_r)}_{n-m}(j;a,b)H^{(s)}_m(xr\log c;\lambda)\right)\frac{t^n}{n!}.
\end{align*}
This completes the proof of the theorem.
\end{proof}
\end{thm}

\smallskip
The next theorem contains an explicit formula for ${B}^{(k_1, k_2,\ldots, k_r)}_n (x; a, b, c)$.

\bigskip
\begin{thm}\label{thm40} {\rm ({\bf Explicit Formula})} For $k \in \mathbb{Z}$, $n\geq0$, we have
\begin{equation}
{B}^{(k_1, k_2,\ldots, k_r)}_n (x; a, b, c) =\sum_{m_r>\ldots >m_1>0}\sum_{j=0}^{m_r-r}\frac{(-1)^j\binom{m_r-r}{j}(rx-j\ln a-(j+1)\ln b)^n}{m_1^{k_1}m_2^{k_2}\ldots m_r^{k_r}}.
\end{equation}
\begin{proof}
\begin{align*}
\frac{{\rm Li}_{(k_1, k_2,\ldots, k_r)}(1-(ab)^{-t})}{(b^t-a^{-t})^r}&=b^{-rt}\left(\sum_{m_r>\ldots >m_1>0}\frac{(1-(ab)^{-t})^{m_r-r}}{m_1^{k_1}m_2^{k_2}\ldots m_r^{k_r}}\right)\\
&=b^{-rt}\sum_{m_r>\ldots >m_1>0}\frac{1}{m_1^{k_1}m_2^{k_2}\ldots m_r^{k_r}}\sum_{j=0}^{m_r-r}(-1)^j\binom{m_r-r}{j}e^{-jt\ln (ab)}\\
&=\sum_{m_r>\ldots >m_1>0}\frac{1}{m_1^{k_1}m_2^{k_2}\ldots m_r^{k_r}}\sum_{j=0}^{m_r-r}(-1)^j\binom{m_r-r}{j}e^{-t(j\ln a+(j+1)\ln b)}.
\end{align*}
So, we get

\smallskip
\begin{equation*}
\frac{{\rm Li}_{(k_1, k_2,\ldots, k_r)}(1-(ab)^{-t})}{(b^t-a^{-t})^r}e^{xrt\ln c}\qquad\qquad\qquad\qquad\qquad\qquad\qquad\qquad\qquad\qquad\qquad\qquad\qquad\qquad\qquad\qquad\qquad\qquad
\end{equation*}
\begin{align*}
&=\sum_{m_r>\ldots >m_1>0}\frac{1}{m_1^{k_1}m_2^{k_2}\ldots m_r^{k_r}}\sum_{j=0}^{m_r-r}(-1)^j\binom{m_r-r}{j}e^{t(rx-j\ln a-(j+1)\ln b)}\\
&=\sum_{n=0}^{\infty}\left(\sum_{m_r>\ldots >m_1>0}\frac{1}{m_1^{k_1}m_2^{k_2}\ldots m_r^{k_r}}\sum_{j=0}^{m_r-r}(-1)^j\binom{m_r-r}{j}(rx-j\ln a-(j+1)\ln b)^n\right)\frac{t^n}{n!}
\end{align*}
By comparing the coefficients of $\frac{t^n}{n!}$ on both sides, the proof is completed.
\end{proof}
\end{thm}

\bigskip
The next theorem contains an expression of ${B}^{(k_1, k_2,\ldots, k_r)}_n(x;a,b,c)$ as polynomial in $x$.

\bigskip
\begin{thm}\label{thm41}
The generalized multi poly-Bernoulli polynomials satisfy the following relation
\begin{equation}\label{eqnnew1}
{B}^{(k_1, k_2,\ldots, k_r)}_n(x;a,b,c)=\sum_{i=0}^n\binom{n}{i}(\ln c)^{n-i}{B}^{(k_1, k_2,\ldots, k_r)}_i(a,b)x^{n-i}
\end{equation}
\begin{proof}
Using (\ref{multipolybern}), we have
\begin{eqnarray*}
\sum_{n=0}^{\infty}{B}^{(k_1, k_2,\ldots, k_r)}_n(x;a,b,c)\frac{t^n}{n!}&=&\frac{{\rm Li}_{(k_1, k_2,\ldots, k_r)}(1-(ab)^{-t})}{(b^t-a^{-t})^r}c^{xt}=e^{xt\ln c}\sum_{n=0}^{\infty}{B}^{(k_1, k_2,\ldots, k_r)}_n(a,b)\frac{t^n}{n!}\\
&=&\sum_{n=0}^{\infty}\sum_{i=0}^n\frac{(xt\ln c)^{n-i}}{(n-i)!}{B}^{(k_1, k_2,\ldots, k_r)}_i(a,b)\frac{t^{i}}{i!}\\
&=&\sum_{n=0}^{\infty}\left(\sum_{i=0}^n\binom{n}{i}(\ln c)^{n-i}{B}^{(k_1, k_2,\ldots, k_r)}_i(a,b)x^{n-i}\right)\frac{t^{n}}{n!}.
\end{eqnarray*}
Comparing the coefficients of $\frac{t^{n}}{n!}$, we obtain the desired result.
\end{proof}
\end{thm}

Note that, when $a=c=e$ and $b=1$, Definition \ref{multipolybernoulli} reduces to
\begin{equation}\label{multipolybernpolynomial}
\frac{{\rm Li}_{(k_1, k_2,\ldots, k_r)}(1-e^{-t})}{(1-e^{-t})^r}e^{rxt}=\sum_{n=0}^{\infty}{B}^{(k_1, k_2,\ldots, k_r)}_n(x)\frac{t^n}{n!}.
\end{equation}
The following theorem gives a relation between ${B}^{(k_1, k_2,\ldots, k_r)}_n(x;a,b,c)$ and ${B}^{(k_1, k_2,\ldots, k_r)}_n(x)$.

\bigskip
\begin{thm}\label{thm42}
The generalized multi poly-Bernoulli polynomials satisfy the following relation
\begin{equation}\label{eqnnew2}
{B}^{(k_1, k_2,\ldots, k_r)}_n(x;a,b,c)=(\ln a+\ln b)^n{B}^{(k_1, k_2,\ldots, k_r)}_n\left(\frac{x\ln c-r\ln b}{\ln a+\ln b}\right)
\end{equation}
\begin{proof}
Using (\ref{multipolybern}), we have
\begin{eqnarray*}
\sum_{n=0}^{\infty}{B}^{(k_1, k_2,\ldots, k_r)}_n(x;a,b,c)\frac{t^n}{n!}&=&\frac{{\rm Li}_{(k_1, k_2,\ldots, k_r)}(1-(ab)^{-t})}{b^{rt}(1-(ab)^{-t})^r}e^{xt\ln c}\\
&=&e^{\frac{x\ln c-r\ln b}{\ln ab}t\ln ab}\frac{{\rm Li}_{(k_1, k_2,\ldots, k_r)}(1-e^{-t\ln ab})}{1+e^{-t\ln ab}}\\
&=&\sum_{n=0}^{\infty}(\ln a+\ln b)^n{B}^{(k_1, k_2,\ldots, k_r)}_n\left(\frac{x\ln c-r\ln b}{\ln a+\ln b}\right)\frac{t^{n}}{n!}.
\end{eqnarray*}
Comparing the coefficients of $\frac{t^{n}}{n!}$, we obtain the desired result.
\end{proof}
\end{thm}

\bigskip
\begin{thm}\label{thm43}
The generalized poly-Bernoulli polynomials satisfy the following relation
\begin{equation}\label{eqnnew3}
\frac{d}{dx}{B}^{(k_1, k_2,\ldots, k_r)}_{n+1}(x;a,b,c)=(n+1)(\ln c){B}^{(k_1, k_2,\ldots, k_r)}_{n}(x;a,b,c)
\end{equation}
\begin{proof}
Using (\ref{multipolybern}), we have
\begin{eqnarray*}
\sum_{n=0}^{\infty}\frac{d}{dx}{B}^{(k_1, k_2,\ldots, k_r)}_{n}(x;a,b,c)\frac{t^{n}}{n!}&=&\frac{t(\ln c){\rm Li}_{(k_1, k_2,\ldots, k_r)}(1-(ab)^{-t})}{(b^t-a^{-t})^r}e^{xrt\ln c}\\
\sum_{n=0}^{\infty}\frac{d}{dx}{B}^{(k_1, k_2,\ldots, k_r)}_{n}(x;a,b,c)\frac{t^{n-1}}{n!}&=&\sum_{n=0}^{\infty}(r\ln c){B}^{(k_1, k_2,\ldots, k_r)}_{n}(x;a,b,c)\frac{t^{n}}{n!}.
\end{eqnarray*}
Hence, 
\begin{eqnarray*}
\sum_{n=0}^{\infty}\frac{1}{n+1}\frac{d}{dx}{B}^{(k_1, k_2,\ldots, k_r)}_{n+1}(x;a,b,c)\frac{t^{n}}{n!}=\sum_{n=0}^{\infty}(\ln c){B}^{(k_1, k_2,\ldots, k_r)}_{n}(x;a,b,c)\frac{t^{n}}{n!}.
\end{eqnarray*}
Comparing the coefficients of $\frac{t^{n}}{n!}$, we obtain the desired result.
\end{proof}
\end{thm}

The following corollary immediately follows from Theorem \ref{thm43} by taking $c=e$. For brevity, let us denote ${B}^{(k_1, k_2,\ldots, k_r)}_{n}(x;a,b,e)$ by ${B}^{(k_1, k_2,\ldots, k_r)}_{n}(x;a,b)$.

\bigskip
\begin{cor}\label{corr}
The generalized poly-Bernoulli polynomials are Appell polynomials in the sense that
\begin{equation}\label{eqnnew4}
\frac{d}{dx}{B}^{(k_1, k_2,\ldots, k_r)}_{n+1}(x;a,b)=(n+1){B}^{(k_1, k_2,\ldots, k_r)}_{n}(x;a,b)
\end{equation}
\end{cor}

\smallskip
Consequently, using the characterization of Appell polynomials \cite{Lee, Shohat, Toscano}, the following addition formula can easily be obtained.  

\bigskip
\begin{cor}\label{cor2}
The generalized poly-Bernoulli polynomials satisfy the following addition formula
\begin{equation}\label{eqnnew5}
{B}^{(k_1, k_2,\ldots, k_r)}_{n}(x+y;a,b)=\sum_{i=0}^n\binom{n}{i}{B}^{(k_1, k_2,\ldots, k_r)}_i(x;a,b)y^{n-i}
\end{equation}
\end{cor}

However, we can derive the addition formula for ${B}^{(k_1, k_2,\ldots, k_r)}_n(x;a,b,c)$ as follows
\begin{eqnarray*}
\sum_{n=0}^{\infty}{B}^{(k_1, k_2,\ldots, k_r)}_n(x+y;a,b,c)\frac{t^n}{n!}&=&\frac{{\rm Li}_{(k_1, k_2,\ldots, k_r)}(1-(ab)^{-t})}{(b^t-a^{-t})^r}c^{(x+y)rt}\\
&=&\frac{{\rm Li}_{(k_1, k_2,\ldots, k_r)}(1-(ab)^{-t})}{(b^t-a^{-t})^r}c^{xrt}c^{yrt}\\
&=&\left(\sum_{n=0}^{\infty}{B}^{(k_1, k_2,\ldots, k_r)}_n(x;a,b,c)\frac{t^n}{n!}\right)\left(\sum_{n=0}^{\infty}(yr\ln c)^n\frac{t^n}{n!}\right)\\
&=&\sum_{n=0}^{\infty}\left(\sum_{i=0}^n\binom{n}{i}(yr\ln c)^{n-i}{B}^{(k_1, k_2,\ldots, k_r)}_i(x;a,b,c)\right)\frac{t^n}{n!}.
\end{eqnarray*}
Comparing the coefficients of $\frac{t^n}{n!}$ yields the following result.

\bigskip
\begin{thm}\label{multithm2}
The generalized poly-Bernoulli polynomials satisfy the following addition formula
\begin{equation*}
{B}^{(k_1, k_2,\ldots, k_r)}_{n}(x+y;a,b,c)=\sum_{i=0}^n\binom{n}{i}(r\ln c)^{n-i}{B}^{(k_1, k_2,\ldots, k_r)}_{i}(x;a,b,c)y^{n-i}.
\end{equation*}
\end{thm}

\begin{rem}
When $x=0$, we have
\begin{equation*}
{B}^{(k_1, k_2,\ldots, k_r)}_{n}(y;a,b,c)=\sum_{i=0}^n\binom{n}{i}(r\ln c)^{n-i}{B}^{(k_1, k_2,\ldots, k_r)}_{i}(a,b,c)y^{n-i}.
\end{equation*}
When $y=(m-1)x$, we further get
\begin{equation*}
{B}^{(k_1, k_2,\ldots, k_r)}_{n}(mx;a,b,c)=\sum_{i=0}^n\binom{n}{i}(r\ln c)^{n-i}{B}^{(k_1, k_2,\ldots, k_r)}_{i}(x;a,b,c)(m-1)^{n-i}x^{n-i},
\end{equation*}
which is a kind of multiplication formula.
\end{rem}

\bigskip
Now let us consider certain symmetrized generalization of ${B}^{(k_1, k_2,\ldots, k_r)}_{n}(x;a,b,c)$.

\begin{defn}\label{Bdefn5}
For $m, n\ge0$, we define
\begin{align*}
\mathcal{C}^{(m)}_{n,r}(x,y;a,b,c)&=\sum_{k_1+k_2+\ldots +k_r=m}\binom{m}{k_1,k_2,\ldots k_r}\frac{{B}^{(-k_1,-k_2,\ldots -k_{r-1})}_n(x;a,b,c)}{(\ln a+\ln b)^n}\times\\
&\;\;\;\times\left(y\ln c-\frac{(r-1)\ln b}{\ln a +\ln b}\right)^{k_r}.
\end{align*}
\end{defn}

The following theorem contains the double generating function for $C^{(m)}_n(x,y;a,b,c)$.

\bigskip
\begin{thm}\label{Bthmm4}
For $n,m\ge0$, we have
\begin{equation}\label{Beqnnnew6}
\sum_{n=0}^{\infty}\sum_{m=0}^{\infty}\mathcal{C}^{(m)}_{n,r}(x,y;a,b,c)\frac{t^n}{n!}\frac{u^m}{m!}=\frac{e^{\left(y\ln c-\frac{(r-1)\ln b}{\ln a +\ln b}\right)u}e^{\left(x\ln c-\frac{(r-1)\ln b}{\ln a +\ln b}\right)t}e^{\binom{r}{2}u+(r-1)t}}{\prod_{i=1}^{r-1}(e^t+e^{iu}-e^{t+iu})}.
\end{equation}
\begin{proof}
$$\sum_{n=0}^{\infty}\sum_{m=0}^{\infty}\mathcal{C}^{(m)}_n(x,y;a,b,c)\frac{t^n}{n!}\frac{u^m}{m!}\qquad\qquad\qquad\qquad\qquad\qquad\qquad\qquad\qquad\qquad\qquad\qquad\qquad\qquad\qquad\qquad$$
\begin{eqnarray*}
&=&\sum_{n=0}^{\infty}\sum_{m=0}^{\infty}\sum_{k_1+k_2+\ldots +k_r=m}\frac{{B}^{(-k_1,-k_2,\ldots -k_{r-1})}_n(x;a,b,c)}{(\ln a+\ln b)^n}\left(y\ln c-\frac{(r-1)\ln b}{\ln a +\ln b}\right)^{k_r}\frac{t^n}{n!}\times\\
&&\;\;\;\;\;\;\times\frac{u^m}{k_1!k_2!\ldots k_r!}\\
&=&\sum_{n=0}^{\infty}\sum_{k_1+k_2+\ldots +k_r\ge0}\frac{{B}^{(-k_1,-k_2,\ldots -k_{r-1})}_n(x;a,b,c)}{(\ln a+\ln b)^n}\left(y\ln c-\frac{(r-1)\ln b}{\ln a +\ln b}\right)^{k_r}\frac{t^n}{n!}\times\\
&&\;\;\;\;\;\;\times\frac{u^{k_1+k_2+\ldots +k_r}}{k_1!k_2!\ldots k_r!}\\
&=&\sum_{n=0}^{\infty}\sum_{k_1+k_2+\ldots +k_{r-1}\ge0}\frac{{B}^{(-k_1,-k_2,\ldots -k_{r-1})}_n(x;a,b,c)}{(\ln a+\ln b)^n}\sum_{k_r\ge0}\left(y\ln c-\frac{(r-1)\ln b}{\ln a +\ln b}\right)^{k_r}\frac{u^{k_r}}{k_r!}\times\\
&&\;\;\;\;\;\;\times\frac{t^n}{n!}\frac{u^{k_1+k_2+\ldots +k_{r-1}}}{k_1!k_2!\ldots k_{r-1}!}\\
&=&e^{\left(y\ln c-\frac{(r-1)\ln a}{\ln a +\ln b}\right)u}\sum_{n=0}^{\infty}\sum_{k_1+k_2+\ldots +k_{r-1}\ge0}\frac{{B}^{(-k_1,-k_2,\ldots -k_{r-1})}_n(x;a,b,c)}{(\ln a+\ln b)^n}\frac{t^n}{n!}\frac{u^{k_1+k_2+\ldots +k_{r-1}}}{k_1!k_2!\ldots k_{r-1}!}
\end{eqnarray*}
Using identity (\ref{eqnnew2}), we obtain
$$\sum_{n=0}^{\infty}\sum_{m=0}^{\infty}\mathcal{C}^{(m)}_{n,r}(x,y;a,b,c)\frac{t^n}{n!}\frac{u^m}{m!}\qquad\qquad\qquad\qquad\qquad\qquad\qquad\qquad\qquad\qquad\qquad\qquad\qquad\qquad\qquad\qquad$$
\begin{eqnarray*}
&=&e^{\left(y\ln c-\frac{(r-1)\ln b}{\ln a +\ln b}\right)u}\sum_{k_1+k_2+\ldots +k_{r-1}\ge0}\sum_{n=0}^{\infty}{B}^{(-k_1,-k_2,\ldots -k_{r-1})}_n\left(x\ln c-\frac{(r-1)\ln b}{\ln a +\ln b}\right)\frac{t^n}{n!}\times\\
&&\;\;\;\;\;\;\times\frac{u^{k_1+k_2+\ldots +k_{r-1}}}{k_1!k_2!\ldots k_{r-1}!}\\
&=&e^{\left(y\ln c-\frac{(r-1)\ln b}{\ln a +\ln b}\right)u}e^{\left(x\ln c-\frac{(r-1)\ln b}{\ln a +\ln b}\right)t}\sum_{k_1+k_2+\ldots +k_{r-1}\ge0}\frac{{\rm Li}_{(-k_1, -k_2,\ldots, -k_{r-1})}(1-e^{-t})}{(1-e^{-t})^{r-1}}\times\\
&&\;\;\;\;\;\;\times\frac{u^{k_1+k_2+\ldots +k_{r-1}}}{k_1!k_2!\ldots k_{r-1}!}
\end{eqnarray*}
\begin{eqnarray*}
&=&\frac{e^{\left(y\ln c-\frac{(r-1)\ln b}{\ln a +\ln b}\right)u}e^{\left(x\ln c-\frac{(r-1)\ln b}{\ln a +\ln b}\right)t}}{(1-e^{-t})^{r-1}}\sum_{0<m_1<m_2<\ldots <m_{r-1}}(1-e^{-t})^{m_{r-1}}e^{u(m_1+m_2+\ldots +m_{r-1})}\\
&=&\frac{e^{y\ln c-\left(\frac{(r-1)\ln b}{\ln a +\ln b}\right)u}e^{\left(x\ln c-\frac{(r-1)\ln b}{\ln a +\ln b}\right)t}e^{\binom{r}{2}u}}{\prod_{i=1}^{r-1}(1-e^{iu}(1-e^{-t}))}\\
&=&\frac{e^{\left(y\ln c-\frac{(r-1)\ln b}{\ln a +\ln b}\right)u}e^{\left(x\ln c-\frac{(r-1)\ln b}{\ln a +\ln b}\right)t}e^{\binom{r}{2}u+(r-1)t}}{\prod_{i=1}^{r-1}(e^t+e^{iu}-e^{t+iu})}.
\end{eqnarray*}
\end{proof}
\end{thm}

\begin{rem}
We observe that duality relation will not work for $\mathcal{C}^{(m)}_{n,r}(x,y;a,b,c)$ when $r\geq3$. However, when $r=2$, the double generating function in Theorem \ref{Bthmm4} yields
$$\sum_{n=0}^{\infty}\sum_{m=0}^{\infty}\mathcal{C}^{(m)}_{n,2}(x,y;a,b,c)\frac{t^n}{n!}\frac{u^m}{m!}=\frac{e^{\left(y\ln c-\frac{\ln b}{\ln a +\ln b}\right)u}e^{\left(x\ln c-\frac{\ln b}{\ln a +\ln b}\right)t}e^{u+t}}{(e^t+e^{u}-e^{t+u})},$$
which implies the following duality relation
$$\mathcal{C}^{(m)}_{n,2}(x,y;a,b,c)=\mathcal{C}^{(m)}_{n,2}(y,x;a,b,c).$$
Furthermore, when $c=e$, we get 
$$\mathcal{C}^{(m)}_{n}(x,y;a,b)=\mathcal{C}^{(m)}_{n}(y,x;a,b),$$
which is exactly the duality relation that appeared in \cite{CorJol5}.
\end{rem}

\section{Hurwitz-Lerch Type Multi Poly-Bernoulli Polynomials}
Consider the case in which $x=1$, $a=e$ and $b=c=1$ for the parameters in Definition \ref{multipolybernoulli}. Then we have
\begin{equation}\label{multipolybernpolynomial-1}
\frac{{\rm Li}_{(k_1, k_2,\ldots, k_r)}(1-e^{-t})}{(1-e^{-t})^r}=\sum_{n=0}^{\infty}{B}^{(k_1, k_2,\ldots, k_r)}_n\frac{t^n}{n!}.
\end{equation}
This can be generalized using the following generalization of Hurwitz-Lerch multiple zeta values 
\begin{equation}\label{HLerch1}
{ \Phi}_{(k_1,k_2,\ldots, k_r)}(z,a)=\sum_{ 0\le m_1\le m_2\le\ldots \le m_r }\frac{z^{m_r}}{(m_1+a-r+1)^{k_1} (m_2+a-r+2)^{k_2}\ldots (m_r+a)^{k_r}}.
\end{equation}
Note that
\begin{align*}
{\rm Li}_{(k_1,k_2,\ldots, k_r)}(z)&=\sum_{ 0< m_1< m_2<\ldots < m_r }\frac{z^{m_r}}{m_1^{k_1} m_2^{k_2}\ldots m_r^{k_r}}\\
&=z^r\sum_{ 0< m_1< m_2<\ldots < m_r }\frac{z^{m_r-r}}{{m_1}^{k_1} {m_2}^{k_2}\ldots {m_r}^{k_r}}\\
&=z^r\sum_{ 0\le m_1\le m_2\le\ldots \le m_r }\frac{z^{m_r}}{{(m_1+1)}^{k_1} {(m_2+2)}^{k_2}\ldots {(m_r+r)}^{k_r}}\\
&=z^r{ \Phi}_{(k_1,k_2,\ldots, k_r)}(z,r)
\end{align*}
Thus, we have
\begin{equation}
\frac{{\rm Li}_{(k_1,k_2,\ldots, k_r)}(z)}{z^r}={ \Phi}_{(k_1,k_2,\ldots, k_r)}(z,r).
\end{equation}
More precisely, one can generalize (\ref{multipolybernpolynomial-1}) as follows
\begin{equation}\label{HpolyBNew}
{ \Phi}_{(k_1,k_2,\ldots, k_r)}(1-e^{-t},a)=\sum_{n=0}^{\infty}B_{n,a}^{(k_1,k_2,\ldots, k_r)}\frac{t^n}{n!}.
\end{equation}
We call $B_{n,a}^{(k_1,k_2,\ldots, k_r)}$ as {\it Hurwitz-Lerch Type Multi Poly-Bernoulli Numbers}. Furthermore, we can define the {\it Hurwitz-Lerch Type Multi Poly-Bernoulli Polynomials} $B_{n,a}^{(k_1,k_2,\ldots, k_r)}(x)$ as follows
\begin{equation}\label{HpolyB1}
{ \Phi}_{(k_1,k_2,\ldots, k_r)}(1-e^{-t},a)e^{rxt}=\sum_{n=0}^{\infty}B_{n,a}^{(k_1,k_2,\ldots, k_r)}(x)\frac{t^n}{n!}
\end{equation}
where $B_{n,a}^{(k_1,k_2,\ldots, k_r)}(0)=B_{n,a}^{(k_1,k_2,\ldots, k_r)}$.

\smallskip
The next theorem contains an explicit formula for $B_{n,a}^{(k_1,k_2,\ldots, k_r)}(x)$ expressed in terms of the $r$-Whitney numbers of the second kind $W_{m,r}(n,k)$ introduced by I. Mez\H{o} in \cite{Mezo} as coefficients of the following exponential generating function
\begin{equation}\label{rbStirling}
\sum_{n=k}^{\infty}W_{m,r}(n,k)\frac{z^n}{n!}=\frac{e^{rz}}{k!}\left(\frac{e^{mz}-1}{m}\right)^k.
\end{equation}
It is worth mentioning that the ordinary Bernoulli numbers $B_n$ defined by
$$\frac{z}{e^z-1}=\sum_{n=0}^{\infty}B_n\frac{z^n}{n!}$$
are expressed in terms $r$-Whitney numbers of the first and second kind in \cite{Mezo} as follows 
$$\binom{n+1}{l}B_{n-l+1}=\frac{n+1}{m^{n-l+1}}\sum_{k=0}^nW_{m,r}(n,k)\frac{w_{m,r}(k+1,l)}{k+1}$$
where $w_{m,r}(n,k)$ are the $r$-Whitney numbers of the first kind defined in \cite{Mezo} as coefficients of the following exponential generating function
$$\sum_{n=k}^{\infty}w_{m,r}(n,k)\frac{z^n}{n!}=(1+mz)^{\frac{-r}{m}}\frac{\ln^k(1+mz)}{m^kk!}.$$
It is important to note that the $r$-Whitney numbers of the second kind $W_{m,r}(n,k)$ are equivalent to the $(r,\beta)$-Stirling numbers $\braced{n}{m}_{\beta,r}$ in \cite{CORCA}. More precisely,
$$W_{\beta,r}(n,m)=\braced{n}{m}_{\beta,r}.$$ 

\bigskip
\begin{thm}\label{res1}
The Hurwitz-Lerch type multi poly-Bernoulli polynomials have the following explicit formula
\begin{equation}
B_{n,a}^{(k_1,k_2,\ldots, k_r)}(x)=\sum_{0\le m_1\le m_2\le\ldots \le m_r\le n}\frac{m_r!W_{-1,xr}(n,m_r)}{(m_1+a-r+1)^{k_1} (m_2+a-r+2)^{k_2}\ldots (m_r+a)^{k_r}}.
\end{equation}
\begin{proof}
Using (\ref{HLerch1}) and (\ref{rbStirling}), we have
$$\sum_{n=0}^{\infty}B_{n,a}^{(k_1,k_2,\ldots, k_r)}(x)\frac{t^n}{n!}={\Phi}_{(k_1,k_2,\ldots, k_r)}(1-e^{-t},a)e^{xrt}\qquad\qquad\qquad\qquad\qquad\qquad\qquad\qquad\qquad\qquad\qquad\qquad\qquad\qquad$$
\begin{align*}
&=\sum_{ 0\le m_1\le m_2\le\ldots \le m_r }\frac{m_r!}{(m_1+a-r+1)^{k_1} (m_2+a-r+2)^{k_2}\ldots (m_r+a)^{k_r}}\frac{e^{xrt}(e^{-t}-1)^{m_r}}{(-1)^{m_r}m_r!}\\
&=\sum_{ 0\le m_1\le m_2\le\ldots \le m_r }\frac{m_r!}{(m_1+a-r+1)^{k_1} (m_2+a-r+2)^{k_2}\ldots (m_r+a)^{k_r}}\sum_{n=m_r}^{\infty}W_{-1,xr}(n,m_r)\frac{t^n}{n!}\\
&=\sum_{n=0}^{\infty}\sum_{ 0\le m_1\le m_2\le\ldots \le m_r\le n }\frac{m_r!W_{-1,xr}(n,m_r)}{(m_1+a-r+1)^{k_1} (m_2+a-r+2)^{k_2}\ldots (m_r+a)^{k_r}}\frac{t^n}{n!}
\end{align*}
Comparing the coefficients of $\frac{t^n}{n!}$ completes the proof of the theorem.
\end{proof}
\end{thm}

Note that (\ref{rbStirling}) implies 
$$W_{-1,0}(n,m_r)=(-1)^{n+m_r}\braced{n}{m_r}.$$ 
Hence, as a direct consequence of Theorem \ref{res1} with $x=0$, we have the following corollary.

\begin{cor} 
The Hurwitz-Lerch type multi poly-Bernoulli numbers equal
\begin{equation}\label{cor1}
B_{n,a}^{(k_1,k_2,\ldots, k_r)}=\sum_{0\le m_1\le m_2\le\ldots \le m_r\le n}\frac{(-1)^{n+m_r}m_r!\braced{n}{m_r}}{(m_1+a-r+1)^{k_1} (m_2+a-r+2)^{k_2}\ldots (m_r+a)^{k_r}}.
\end{equation}
\end{cor}

\begin{rem}
When $r=1$, equation (\ref{cor1}) gives
$$B_{n,a}^{(k_1)}=\sum_{0\le m_1\le n}\frac{(-1)^{n+m_1}m_1!\braced{n}{m_1}}{(m_1+a)^{k_1}}$$
which is exactly the explicit formula for Hurwitz-Lerch type poly-Bernoulli numbers in \eqref{HLExplicit}.
\end{rem}
\begin{rem}
Note that, when $k_1, k_2, \ldots, k_r$ are all positive integers, $B_{n,a}^{(-k_1,-k_2,\ldots, -k_r)}$ are always nonnegative integers. This implies that we can possibly draw combinatorial interpretations of these numbers.
\end{rem}

\bigskip
\begin{flushleft}
{\bf Roberto B. Corcino}\\
Cebu Normal University\\
Cebu City, Philippines\\
e-mail: rcorcino@yahoo.com
\end{flushleft}

\begin{flushleft}
{\bf Hassan Jolany}\\
Universit\'e des Sciences et Technologies de Lille\\
UFR de Math\'ematiques\\
Laboratoire Paul Painlev\'e\\
CNRS-UMR 8524 59655 Villeneuve d'Ascq Cedex/France\\
e-mail: hassan.jolany@math.univ-lille1.fr
\end{flushleft}

\begin{flushleft}
{\bf Cristina B. Corcino}\\
Cebu Normal University\\
Cebu City, Philippines\\
e-mail: cristinacorcino@yahoo.com
\end{flushleft}

\begin{flushleft}
{\bf Takao Komatsu}\\
School of Mathematics and Statistics\\
Wuhan University\\ 
Wuhan 430072 China\\
e-mail: komatsu@whu.edu.cn
\end{flushleft}

\end{document}